\documentclass[12pt,a4paper]{article}
\usepackage[latin2]{inputenc}
\usepackage{amsmath}
\usepackage{amsfonts}
\usepackage{amssymb}
\usepackage{graphicx}
\usepackage[left=2.00cm, right=2.00cm, top=2.50cm, bottom=2.50cm]{geometry}
\usepackage[T1]{fontenc}

\def\Q{{\mathbb Q}}
\def\Z{{\mathbb Z}}

\newtheorem{lemma}{Lemma}
\newtheorem{theorem}[lemma]{Theorem}

\title{
Totally real Thue inequalities \\ over imaginary quadratic fields 
}
\author{
Istv\'{a}n Ga\'{a}l\thanks{
        Research supported in part by K115479 from the
        Hungarian National Foundation for Scientific Research				
				and by the EFOP-3.6.1-16-2016-00022 project. 
				The project is co-financed by the European Union and the European Social Fund.	
                         },\; \\
{\small University of Debrecen, Mathematical Institute} \\
{\small H--4002 Debrecen Pf.400., Hungary,} 
{\small e--mail: gaal.istvan@unideb.hu},
 \\ \\
Borka Jadrijevi\'c\thanks{
Research supported in part
 by the Croatian Science Foundation under the
project no. 6422.   }\\
{\small University of Split,}
{\small Faculty of Science, }\\
{\small Ru\dj era Bo\v{s}kovi\'{c}a 33, 21000 Split, Croatia,}
{\small e--mail: borka@pmfst.hr}
 \\ \\ 
L\'aszl\'o Remete\thanks{
        Research supported by the \'UNKP-17-3 new national excellence program of the Ministry of human capacities.}\; \\
{\small University of Debrecen, Mathematical Institute} \\
{\small H--4002 Debrecen Pf.400., Hungary,} 
{\small e--mail: remete.laszlo@science.unideb.hu}
}

\begin{document}
\baselineskip=17pt

\maketitle
\thispagestyle{empty}

\renewcommand{\thefootnote}{\arabic{footnote}}
\setcounter{footnote}{0}


\begin{abstract}
Let $F(x,y)$ be an irreducible binary form of degree $\geq 3$ with
integer coefficients and with real roots. Let $M$ be an imaginary quadratic field,
with ring of integers $\Z_M$. Let $K>0$.
We describe an efficient method how to reduce the resolution of the relative Thue 
inequalities 
\[
|F(x,y)|\leq K  \;\; (x,y\in \Z_M)
\]
to the resolution of absolute Thue inequalities of type
\[
|F(x,y)|\leq k \;\; (x,y\in \Z).
\]
We illustrate our method with an explicit example.
\end{abstract}

\newpage

\section{Introduction}

Let $F(x,y)\in\Z[x,y]$ be an irreducible binary form of degree $\geq 3$
and let $a\in\Z\setminus \{0\}$.
There is an extensive literature of {\bf Thue equations} of type 
\[
F(x,y)=a\;\;{\rm in}\;\;x,y\in\Z.
\]
In 1909 A. Thue \cite{thue} proved that these equations admit only finitely many solutions.
In 1967 A. Baker \cite{baker} gave effective upper bounds for the solutions.
Later on authors constructed numerical methods to reduce the bounds and 
to explicitly calculate the solutions, see \cite{book} for a summary.

Let $M$ be an algebraic number field with ring of integers $\Z_M$.
Let $F(x,y)\in\Z_M[x,y]$ be  an irreducible binary form of degree 
$n\geq 3$ and let $\mu\in\Z_M\setminus \{0\}$. As a generalization of Thue equations
consider {\bf relative Thue equations} of type
\[
F(x,y)=\mu\;\;{\rm in}\;\;x,y\in\Z_M.
\]

Using Baker's method S. V. Kotov and V. G. Sprindzuk \cite{ks} were first 
to give effective upper
bounds for the solutions of relative Thue equations.
Their theorem has been extended by several authors.
Applying Baker's method, reduction and enumeration algorithms I. Ga\'al and M. Pohst \cite{gp}
gave an efficient algorithm for solving relative Thue equations (see also \cite{book}).

\vspace{0.5cm}

Let $M$ be an imaginary quadratic number field.
Assuming in addition that the roots of $F(x,1)$ are all real,
in the present paper we give an efficient algorithm to reduce
the resolution of {\bf relative Thue inequalities} of the type
\[
|F(x,y)|\leq K \;\;{\rm in}\;\;x,y\in\Z_M
\]
to the resolution of (absolute) {\bf Thue inequalities} of the type 
\[
|F(x,y)|\leq k \;\;{\rm in}\;\;x,y\in\Z.
\]
To find the solutions of the above absolute Thue inequality 
one can use Kash \cite{kash} or Magma \cite{magma} which admit
efficient algorithms for solving (absolute) Thue equations
$F(x,y)=k'$ for $k'\in\Z$ with $|k'|\leq k$.
For an efficient method for calculating "small" solutions of 
Thue inequalities we refer to \cite{thueineq}.

Our method is illustrated with an explicit example.

\section{The main result}

Let $F(x,y)$ be a binary form of degree $n\geq 3$ with rational integer coefficients.
Assume that $f(x)=F(x,1)$ has leading coefficient 1 and distinct real roots $\alpha_1,\ldots,\alpha_n$. 
Let $0<\varepsilon<1,\;\; 0<\eta<1$ and let $K\geq 1$. Set

\begin{eqnarray*}
&A={\displaystyle\min_{i\not=j}|\alpha_i-\alpha_j|},\;\; 
&B={\displaystyle\min_i\prod_{j\not=i} |\alpha_j-\alpha_i|},
\\
&C={\displaystyle \max\left\{\frac{K}{(1-\varepsilon)^{n-1}B},1\right\}},&
\\
&C_1={\displaystyle\max\left\{\frac{K^{1/n}}{\varepsilon A},\;
(2C)^{1/(n-2)}\right\}},\;\;
&C_2={\displaystyle\max\left\{\frac{K^{1/n}}{\varepsilon A},\; 
C^{1/(n-2)}\right\}},
\\
&D={\displaystyle\left(\frac{K}{\eta(1-\varepsilon)^{n-1}AB}\right)^{1/n}},\;
&E={\displaystyle\frac{(1+\eta)^{n-1}K}{(1-\varepsilon)^{n-1}}}.
\end{eqnarray*}

\vspace{1cm}

\noindent
Let $m>1$ be a squarefree positive integer, and set $M=\Q(i\sqrt{m})$. 
Consider the relative Thue inequality
\begin{equation}
|F(x,y)|\leq K \;\; {\rm in}\;\; x,y\in\Z_M.
\label{1}
\end{equation}

\vspace{1cm}

\noindent
If $m\equiv 3\; (\bmod \; 4)$, then $x,y\in\Z_M$ can be written as
\[
x=x_1+x_2\frac{1+i\sqrt{m}}{2}=\frac{(2x_1+x_2)+x_2i\sqrt{m}}{2}, 
\]
\[
y=y_1+y_2\frac{1+i\sqrt{m}}{2}=\frac{(2y_1+y_2)+y_2i\sqrt{m}}{2}
\]
with $x_1,x_2,y_1,y_2\in\Z$.

\noindent
If $m\equiv 1,2\; (\bmod \; 4)$, then $x,y\in\Z_M$ can be written as
\[
x=x_1+x_2i\sqrt{m},\;\; y=y_1+y_2i\sqrt{m}
\]
with $x_1,x_2,y_1,y_2\in\Z$.

\vspace{1cm}

\begin{theorem}
\label{th1}
Let $(x,y)\in\Z_M^2$ be a solution of (\ref{1}). 
Assume that 
\begin{eqnarray}
|y|> C_1 & \;\; {\rm if} \;\; &m\equiv 3\; (\bmod \; 4),\label{f1} \\
|y|> C_2 & \;\; {\rm if} \;\; &m\equiv 1,2\; (\bmod \; 4). \label{f2}
\end{eqnarray}
Then 
\begin{equation}
x_2y_1=x_1y_2.
\label{xy}
\end{equation}

\vspace{0.5cm}
\noindent
I. Further, if $m\equiv 3\; (\bmod \; 4)$, then the following holds:

\vspace{0.5cm}
\noindent
IA1. If  $2y_1+y_2=0$, then  $2x_1+x_2=0$ and
\begin{equation}
|F(x_2,y_2)|\leq \displaystyle{\frac{2^nK}{(\sqrt{m})^n}}.
\label{IA1}
\end{equation}
IA2. If  $|2y_1+y_2|\geq 2D$,  then 
\begin{equation}
|F(2x_1+x_2,2y_1+y_2)|\leq 2^n E.
\label{IA2}\\
\end{equation}
IB1. If $y_2=0$ then $x_2=0$ and
\begin{equation}
|F(x_1,y_1)|\leq K.
\label{IB1}
\end{equation}
IB2. If $|y_2|\geq \displaystyle{\frac{2}{\sqrt{m}} D}$, then
\begin{equation}
|F(x_2,y_2)|\leq \frac{2^n}{(\sqrt{m})^n} E.
\label{IB2}
\end{equation}

\vspace{0.8cm}
\noindent
II. If $m\equiv 1,2\; (\bmod \; 4)$, then the following holds:

\vspace{0.5cm}
\noindent
IIA1. If $y_1=0$ then $x_1=0$ and
\begin{equation}
|F(x_2,y_2)|\leq \displaystyle{\frac{K}{(\sqrt{m})^n}}.
\label{IIA1}
\end{equation}
IIA2.  If   $|y_1|\geq D$,  then
\begin{equation}
|F(x_1,y_1)|\leq E.
\label{IIA2}
\end{equation}
IIB1. If $y_2=0$ then $x_2=0$ and
\begin{equation}
|F(x_1,y_1)|\leq K.
\label{IIB1}
\end{equation}
IIB2. If $|y_2|\geq \displaystyle{\frac{D}{\sqrt{m}}}$, then
\begin{equation}
|F(x_2,y_2)|\leq \frac{E}{(\sqrt{m})^n}.
\label{IIB2}
\end{equation}
\end{theorem}

\vspace{0.5cm}

\noindent
Our result is a far reaching generalization of an idea of \cite{complex6}.

\vspace{1cm}

\section{Proof of the main result}

In the proof of Theorem \ref{th1} we shall use the following Lemma.

\begin{lemma}
\label{lemma1}
Let $x,y\in\Z,y\not=0$. Assume that 
\[
\left|\alpha_{i_0}-\frac{x}{y}\right|\leq \frac{d}{|y|^n}
\]
for some $i_0\; (1\leq i_0\leq n)$ and $d>0$. If
\[
 |y|\geq \left(\frac{d}{\eta A}  \right)^{1/n},
\]
then
\[
\;\; 
|F(x,y)|\leq d(1+\eta)^{n-1}\prod_{j\not={i_0}} |\alpha_j-\alpha_{i_0}|.
\]
\end{lemma}

\vspace{1cm}

\noindent
{\bf Proof of Lemma \ref{lemma1}}\\
By our assumption, we have
\[
\left|\alpha_j-\frac{x}{y}\right|\leq |\alpha_j-\alpha_{i_0}|+\left|\alpha_{i_0}-\frac{x}{y}\right|
\leq (1+\eta) |\alpha_j-\alpha_{i_0}|
\]
for $j\not={i_0}$. Therefore
\[
\prod_{j=1}^n \left|\alpha_j-\frac{x}{y}\right|= 
\left|\alpha_{i_0}-\frac{x}{y}\right| \cdot \prod_{j\not ={i_0}}^n \left|\alpha_j-\frac{x}{y}\right|\leq 
\frac{d}{|y|^n}\cdot (1+\eta)^{n-1}\cdot \prod_{j\not={i_0}} |\alpha_j-\alpha_{i_0}|,
\]
which implies our assertion. \hfill $\Box$

\vspace{1cm}

\noindent
{\bf Proof of Theorem \ref{th1}.}\\

Let $(x,y)\in\Z_M^2$ be an arbitrary solution of (\ref{1}) with $y\neq 0$.
Let $\beta_j=x-\alpha_j y,\; j=1,\ldots,n$, then the inequality (\ref{1}) can be written as
\begin{equation}
|\beta_1\cdots\beta_n|\leq K.
\label{e}
\end{equation}
Let ${i_0}$ be the index with
\[
|\beta_{i_0}|=\min_j |\beta_j|.
\]
Then $|\beta_{i_0}|\leq K^{\frac{1}{n}}$ and together with
 (\ref{f1}) and (\ref{f2}) we get
\[
|\beta_j|\geq |\beta_j-\beta_{i_0}|-|\beta_{i_0}|\geq |\alpha_j-\alpha_{i_0}|\cdot |y|-K^{\frac{1}{n}}
\geq (1-\varepsilon) \cdot |\alpha_j-\alpha_{i_0}|\cdot |y|
\]
for $j\not={i_0}$.
From the previous inequality and (\ref{e}), we have
\begin{equation}
|\beta_{i_0}|\leq \frac{K}{\prod_{j\not={i_0}}|\beta_j|}\leq \frac{c}{|y|^{n-1}}
\label{betai}
\end{equation}
with
\[
c=\frac{K}{(1-\varepsilon)^{n-1}\prod_{j\not={i_0}}|\alpha_j-\alpha_{i_0}|}.
\]
By (\ref{betai}) we obtain
\[
\left|\alpha_{i_0}-\frac{x\overline{y}}{|y|^2}\right|=\left|\alpha_{i_0}-\frac{x}{y}\right|\leq \frac{c}{|y|^{n}},
\]
hence
\[
\left|\alpha_{i_0}|y|^2-x\overline{y}\right|\leq \frac{c}{|y|^{n-2}},
\]
which implies 
\[
|{\rm Im} (x\overline{y})|\leq \frac{c}{|y|^{n-2}}.
\]
Note that $\frac{c}{|y|^{n-2}}<\frac{1}{2}$ and $\frac{c}{|y|^{n-2}}<1$ for 
$m\equiv 3\; (\bmod \; 4)$ and $m\equiv 1,2\; (\bmod \; 4)$,
respectively, according to (\ref{f1}) and (\ref{f2}).
Therefore 
$|{\rm Im} (x\overline{y})|=\frac{1}{2}|x_2y_1-x_1y_2|\sqrt{m}<\frac{1}{2}$
and \\
$|{\rm Im} (x\overline{y})|=|x_2y_1-x_1y_2|\sqrt{m}<1$
for 
$m\equiv 3\; (\bmod \; 4)$ and $m\equiv 1,2\; (\bmod \; 4)$, respectively.
Hence in both cases we have (\ref{xy}).

\vspace{0.5cm}

\noindent
I. Let $m\equiv 3\; (\bmod \; 4)$. 

\noindent
IA. The inequality (\ref{betai}) implies 
$|{\rm Re}(\beta_{i_0})|\leq \frac{c}{|y|^{n-1}}$, i.e.
\begin{equation}
|(2x_1+x_2)-\alpha_{i_0}(2y_1+y_2)|\leq \frac{2c}{|y|^{n-1}}.
\label{h1}
\end{equation}

\noindent
IA1. If $2y_1+y_2=0$, then (\ref{h1}) yields $2x_1+x_2=0$,
and the inequality (\ref{1}) has the form
\[
\left|F\left(\frac{x_2i\sqrt{m}}{2},\frac{y_2i\sqrt{m}}{2}\right)\right|\leq K
\]
whence we get (\ref{IA1}).

\noindent
IA2. If $2y_1+y_2\not=0$, then 
\[
|(2x_1+x_2)-\alpha_{i_0}(2y_1+y_2)|\leq \frac{2c}{|y|^{n-1}}
=\frac{2c}{\displaystyle{\left|\frac{(2y_1+y_2)+y_2i\sqrt{m}}{2}\right|^{n-1}}}
\leq \frac{2^nc}{|2y_1+y_2|^{n-1}}.
\]
Since we have assumed 
\[
|2y_1+y_2|\geq \left(\frac{2^nc}{\eta A}\right)^{1/n},
\]
Lemma \ref{lemma1} implies
\[
|F(2x_1+x_2,2y_1+y_2)|\leq 2^nc(1+\eta)^{n-1}\prod_{j\not={i_0}}|\alpha_j-\alpha_{i_0}|
\]
whence we get (\ref{IA2}).

\noindent
IB. By the inequality (\ref{betai}), we have 
$|{\rm Im}(\beta_{i_0})|\leq \frac{c}{|y|^{n-1}}$, i.e.
\begin{equation}
\sqrt{m}|x_2-\alpha_{i_0}y_2|\leq \frac{2c}{|y|^{n-1}}.
\label{h2}
\end{equation}

\noindent
IB1. If $y_2=0$, then (\ref{h2}) implies $x_2=0$
and the inequality (\ref{1}) has the form
\[
\left|F\left(\frac{2x_1}{2},\frac{2y_1}{2}\right)\right|\leq K
\]
whence we get (\ref{IB1}).

\noindent
IB2. If $y_2\not=0$, then 
\[
|x_2-\alpha_{i_0}y_2|\leq \frac{2c}{\sqrt{m}|y|^{n-1}}
=\frac{2c}{\sqrt{m}\displaystyle{\left|\frac{(2y_1+y_2)+y_2i\sqrt{m}}{2}\right|^{n-1}}}
\leq \frac{2^nc}{(\sqrt{m})^{n}|y_2|^{n-1}}.
\]
Since 
\[
|y_2|\geq \left(\frac{2^nc}{(\sqrt{m})^n\eta A}\right)^{1/n},
\]
Lemma \ref{lemma1} implies
\[
|F(x_2,y_2)|\leq
\frac{2^nc}{(\sqrt{m})^n}(1+\eta)^{n-1}\prod_{j\not={i_0}}|\alpha_j-\alpha_{i_0}|
\]
which implies (\ref{IB2}).

\vspace{1cm}

\noindent
II. Let $m\equiv 1,2\; (\bmod \; 4)$. 

\noindent
IIA.  The inequality (\ref{betai}) implies 
$|{\rm Re}(\beta_{i_0})|\leq \frac{c}{|y|^{n-1}}$, i.e.
\begin{equation}
|x_1-\alpha_{i_0}y_1|\leq \frac{c}{|y|^{n-1}}.
\label{k1}
\end{equation}

\noindent
IIA1. If $y_1=0$, then (\ref{k1}) yields $x_1=0$ and 
the inequality (\ref{1}) has the form
\[
|F(i\sqrt{m}x_2,i\sqrt{m}y_2)|\leq K,
\]
whence we get (\ref{IIA1}).

\noindent
IIA2. If $y_1\not=0$, then 
\[
|x_1-\alpha_{i_0}y_1|\leq \frac{c}{|y|^{n-1}}
=\frac{c}{|y_1+i\sqrt{m}y_2|^{n-1}}
\leq \frac{c}{|y_1|^{n-1}}.
\]
Since we have assumed
\[
|y_1|\geq \left(\frac{c}{\eta A}\right)^{1/n},
\]
Lemma \ref{lemma1} implies
\[
|F(x_1,y_1)|\leq c(1+\eta)^{n-1}\prod_{j\not={i_0}}|\alpha_j-\alpha_{i_0}|
\]
whence we get (\ref{IIA2}).

\noindent
IIB. By the inequality (\ref{betai}) we have
$|{\rm Im}(\beta_{i_0})|\leq \frac{c}{|y|^{n-1}}$, i.e.
\begin{equation}
\sqrt{m}|x_2-\alpha_{i_0}y_2|\leq \frac{c}{|y|^{n-1}}.
\label{k2}
\end{equation}

\noindent
IIB1. If $y_2=0$, then (\ref{k2}) implies  $x_2=0$ 
and the inequality (\ref{1}) has the form
\[
|F(x_1,y_1)|\leq K
\]
which is just our assertion (\ref{IIB1}).

\noindent
IIB2. If $y_2\not=0$, then 
\[
|x_2-\alpha_{i_0}y_2|\leq \frac{c}{\sqrt{m}|y|^{n-1}}
=\frac{c}{|y_1+i\sqrt{m}y_2|^{n-1}}
\leq \frac{c}{(\sqrt{m})^{n}|y_2|^{n-1}}.
\]
Since
\[
|y_2|\geq \left(\frac{c}{(\sqrt{m})^n\eta A}\right)^{1/n},
\]
Lemma \ref{lemma1} implies
\[
|F(x_2,y_2)|
\leq\frac{c}{(\sqrt{m})^n}(1+\eta)^{n-1}\prod_{j\not={i_0}}|\alpha_j-\alpha_{i_0}|
\]
whence we get (\ref{IIB2}).
\hfill $\Box$

\section{How to apply Theorem \ref{th1}}

In this section we give useful hints for a practical application of Theorem \ref{th1}.\\

\noindent
Using the same notation let us consider again 
the relative Thue inequality (\ref{1}).
We describe our algorithm in the case I (for $m\equiv 3 \; (\bmod \; 4)$)
since the case II is completely similar.

\begin{enumerate}
\item
If $|y|\leq C_1$ then we have only finitely many possible
values for $y$ and hence for $y_1,y_2$, as well. 
For each possible $y$ and for all integers $\mu\in\Z_M$ with $|\mu|\leq K$
we calculate the roots of the equation $F(x,y)-\mu=0$ in $x$.
For such a root $x$ we calculate the corresponding  $x_1,x_2$. 
If $x_1,x_2$ are integers, then 
$x\in\Z_M$ and $(x,y)$ is a solutions of (\ref{1}). \\
Alternatively, by 
$|\beta_{i_0}|\leq K^{\frac{1}{n}}$ we obtain $|x|\leq K^{\frac{1}{n}}+
\max|\alpha_j|\cdot C_1$. We can simply enumerate and test the finitely many
possible values of $x_1,x_2$ and $y_1,y_2$.
\item
Assume that $|y|>C_1$. 
  \begin{enumerate}
  \item
If $|2y_1+y_2|<2D$, then 
     \begin{enumerate}
      \item
      If $|y_2|<2D/\sqrt{m}$, then we have only finitely many values for $y_1,y_2$,
			we proceed as in 1.
      \item 
			If $|y_2|\geq 2D/\sqrt{m}$, then we use IB2. We solve 
      $F(x_2,y_2)=k$ for all $k\in\Z$ with $|k|\leq 2^nE/(\sqrt{m})^n$. We determine the possible
			values of $y_1$ which satisfy $|2y_1+y_2|<2D$. 
			We substitute $x_2,y_1,y_2$ into $x_2y_1=x_1y_2$
			to see if there exist corresponding integer $x_1$.
     \end{enumerate}
  \item
	If $|2y_1+y_2|\geq 2D$, then we use IA2. We calculate the solutions
	$X=2x_1+x_2,Y=2y_1+y_2$ of 
      $F(X,Y)=k$ for all $k\in\Z$ with $|k|\leq 2^nE$.
	  \begin{enumerate}
      \item	
			If $|y_2|<2D/\sqrt{m}$ then there are only finitely many possible values for $y_2$.
			We determine $y_1$ from $Y$.
			Using $X=2x_1+x_2$ we set $x_2=X-2x_1$, substitute $x_2=X-2x_1,y_1,y_2$
			into $x_2y_1=x_1y_2$ and test if there is a corresponding $x_1$ in $\Z$.
		 \item	
			If $|y_2|\geq 2D/\sqrt{m}$ we use IB2. We solve 
			$F(x_2,y_2)=k$ for $|k|\leq 2^nE/(\sqrt{m})^n$. 
			We determine $x_1,y_1$ from $x_2,y_2$ and $X,Y$.
     \end{enumerate}
\end{enumerate}
\end{enumerate}

For solving absolute Thue equations $F(x,y)=k$ for certain values $k\in\Z$
one can efficiently apply Kash \cite{kash} and Magma \cite{magma}.

We remark that an appropriate choice of the parameters $\varepsilon,\eta$
of Thereom \ref{th1}
makes the resolution much easier. 
It is worthy to keep $C_1,C_2$ and also $D$ small, 
to avoid extensive tests of small
possible solutions. On the other hand, if $E$ is small, then 
there are fewer Thue equations (over $\Z$) to be solved.
Of course we can not make all these 
constants simultaneously small, therefore we need to make a compromise,
taking into consideration also 
the value of $K$ (which also determines the number of 
Thue equations to be solved). Usually it is worthy to try several
values of $\varepsilon,\eta$ before we start solving (\ref{1}).

\section{An example}

Let $M=\Q(i\sqrt{5})$, and let
\[
F(x,y)=x^4-9x^3y-21x^2y^2+88xy^3+48y^4
\]
and consider the solutions of
\begin{equation}
|F(x,y)|\leq 20 \;\; {\rm in}\;\; x,y\in\Z_M.
\label{pl}
\end{equation}
The polynomial $F(x,y)$ is irreducible and 
the roots of $F(x,1)$ are approximately
\[
-3.4271, -0.49938, 2.7581, 10.1684.
\]
We may set $A=2.9278$, $B=101.7426$. Further, let $\varepsilon=0.1$ and $\eta=0.1$.
We are in case II. Calculating the constants, Theorem \ref{th1} gives:\\
\\
Assume $|y|> 7.2229$. Then:
\begin{eqnarray*}
{\rm IIA1.}&\; {\rm If} \;y_1=0, \;{\rm then}\; x_1=0 \;{\rm and} \;
|F(x_2,y_2)|\leq 0.8000.
\\
{\rm IIA2.}&\; {\rm If} \;|y_1|\geq 0.9796,\; {\rm then}\;
|F(x_1,y_1)|\leq 36.5157.
\\
{\rm IIB1.}&\; {\rm If} \; y_2=0,\;{\rm  then}\; x_2=0\;{\rm  and}\;
|F(x_1,y_1)|\leq 20.
\\
{\rm IIB2.}&\; {\rm If}\; |y_2|\geq 0.4381,\; {\rm then}\;
|F(x_2,y_2)|\leq 1.4606.
\end{eqnarray*}\\

First we consider the values with $|y|\leq C_2=7.2229$.
We have $|x|\leq 20^{\frac{1}{4}}+\max|\alpha_j|\cdot C_2=75.64$.
Enumerating and testing all possible $x=x_1+i\sqrt{5}x_2$ and $y=y_1+i\sqrt{5}y_2$
satisfying these bounds we obtain the solutions 
$(x_1,x_2,y_1,y_2)=(0,0,0,0),(1,0,0,0),(2,0,0,0),(1,0,-2,0)$, $(2,0,-4,0)$,
up to sign.

If $y_1=0$ then by IIA1 we have $x_1=0$ and $|F(x_2,y_2)|\leq 0.8$, whence 
$|F(x_2,y_2)|=0$, $x_2=0,y_2=0$.

If $y_2=0$ then by IIB1 we have $x_2=0$ and $|F(x_1,y_1)|\leq 20$. 
Using Magma we solve $F(x_1,y_1)=k$ 
for $-20\leq k\leq 20$. We obtain the solutions 
$(x_1,y_1)=(0,0),(1,0),(1,-2),(2,0)$, $(2,-4)$, up to sign. 
These bring the above solutions $(x_1,x_2,y_1,y_2)$ again.

From now on we assume that $y_1\neq 0$ and $y_2\neq 0$.

If $|y_1|\leq 0.9796$ and $|y_2|\leq 0.4381$ then by IIA2 we have
$|F(x_1,y_1)|\leq 36.5157$ and by IIB2 we have 
$|F(x_2,y_2)|\leq 1.4606$.
In addition to the above calculation we solve $F(x_1,y_1)=k$
for $21\leq |k|\leq 36$ but we do not get any further solutions. 
Hence the solutions of $|F(x_1,y_1)|\leq 36.5157$ are 
$(x_1,y_1)=(0,0),(1,0),(1,-2),(2,0),(2,-4)$, up to sign.
Also the solutions of $|F(x_2,y_2)|\leq 1.4606$ are
$(x_1,y_1)=(0,0),(1,0),(1,-2)$, up to sign. 
Testing these possible 
$(x_1,x_2,y_1,y_2)$ we do not get any new solutions.

If either $|y_1|< 0.9796$ or $|y_2|< 0.4381$ then $y_1=0$ or $y_2=0$
which cases we have already considered.

Hence all solutions of (\ref{pl}) are
$(x,y)=(0,0),(1,0),(2,0),(1,-2),(2,-4)$,
up to sign. 
The calculation takes just a few seconds.

\end{document}